\documentclass[brochure,12pt,french]{smfbourbaki}
        
\usepackage[T1]{fontenc}
\usepackage{lmodern,amssymb,bm}

\usepackage{babel}

\usepackage[utf8]{inputenc}

\usepackage[colorlinks=true, linkcolor=blue, citecolor=red,
urlcolor=blue]{hyperref}

\usepackage[
backend=bibtex,
style=authoryear, 
citestyle=authoryear-comp,
maxnames=7,
sortcites=false 
]{biblatex}
\usepackage{csquotes}

\usepackage{pifont}
\usepackage{amscd}

\newcommand{\pt}{{\tiny{\ding{71}}}}
\newcommand{\Gt}{\ding{71}}

\newcommand{\R}{\mathbb{R}}
\newcommand{\C}{\mathbb{C}}
\newcommand{\G}{\mathbb{G}}
\newcommand{\N}{\mathbb{N}}
\newcommand{\Z}{\mathbb{Z}}
\newcommand{\T}{\mathbb{T}}

\newcommand{\gA}{\mathfrak{A}}

\newcommand{\cA}{\mathcal{A}}

\newcommand{\cF}{\mathcal{F}}

\newcommand{\cG}{\mathcal{G}}
\newcommand{\cH}{\mathcal{H}}
\newcommand{\cI}{\mathcal{I}}
\newcommand{\cJ}{\mathcal{J}}
\newcommand{\cK}{\mathcal{K}}
\newcommand{\cL}{\mathcal{L}}
\newcommand{\cM}{\mathcal{M}}
\newcommand{\cN}{\mathcal{N}}
\newcommand{\cO}{\mathcal{O}}

\newcommand{\cR}{\mathcal{R}}
\newcommand{\cS}{\mathcal{S}}
\newcommand{\cT}{\mathcal{T}}
\newcommand{\cU}{\mathcal{U}}

\newcommand{\cX}{\mathcal{X}}

\newtoggle{tnbcbx@howcited}
\DeclareEntryOption[boolean]{howcited}[true]{\settoggle{tnbcbx@howcited}{#1}}
\DeclareBibliographyOption[boolean]{howcited}[true]{\settoggle{tnbcbx@howcited}{#1}}
\DeclareTypeOption[boolean]{howcited}[true]{\settoggle{tnbcbx@howcited}{#1}}

\newbibmacro{howcited}{%
  \iftoggle{tnbcbx@howcited}
    {\iffieldundef{shorthand}
       {}
       {\setunit{\addspace}%
        \printtext[parens]{%
          \bibstring{citedas}%
          \setunit{\addcolon\space}%
          \printfield{shorthand}%
          \setunit{\addslash}%
          }}}
    {}}

\renewbibmacro{finentry}{\usebibmacro{howcited}\finentry}

\DefineBibliographyStrings{french}{citedas    = {cité ci-dessus avec l'acronyme}}
\DefineBibliographyStrings{english}{citedas    = {heretofore cited as}}

\DefineBibliographyExtras{french}{\restorecommand\mkbibnamefamily}

\DeclareDelimFormat{nameyeardelim}{\addcomma\space}

\DeclareNameAlias{sortname}{given-family}

\renewcommand{\bibnamedash}{\leavevmode\raise3pt\hbox to3em{\hrulefill}\space}

\AtEveryBibitem{%
  \clearfield{issn} 
  \clearfield{isbn} 
  \clearfield{doi} 
  \clearlist{language} 
  \ifentrytype{online}{}{
  \ifentrytype{unpublished}{}{
    \clearfield{url}
  }
  }
}

\renewbibmacro{in:}{%
    \ifentrytype{article}{}{\printtext{\bibstring{in}\intitlepunct}}}

\DeclareFieldFormat[article,periodical,inreference]{number}{\mkbibparens{#1}}
\DeclareFieldFormat[article,periodical,inreference]{volume}{\mkbibbold{#1}}
\renewbibmacro*{volume+number+eid}{%
    \printfield{volume}%
    \setunit*{\addthinspace}
    \printfield{number}%
    \setunit{\addcomma\space}%
    \printfield{eid}}

\DeclareFieldFormat[article,inbook,incollection]{title}{\enquote{#1}\addcomma} 

\date{Mars 2025}
\bbkannee{77\textsuperscript{e} année, 2024--2025}  
\bbknumero{1236}                                      

\title{Hypoellipticité de polynômes de champs de vecteurs et conjectures de Helffer et Nourrigat}
\subtitle{d'après I. Androulidakis, O. Mohsen, R. Yuncken}

\author{Claire Debord}
\address{IMJ-PRG,  Université Paris Cité \\ 
Bâtiment Sophie Germain \\
8 place Aurélie Nemours \\
Boite Courrier 7012 \\
75205 Paris Cedex 13}
\email{claire.debord@imj-prg.fr}

\begin{document}

\maketitle


\section*{Introduction}

Un opérateur différentiel $D$ sur une variété $M$ est dit \emph{hypoelliptique} si l'équation différentielle $Du=f$ satisfait la propriété de régularité suivante : pour tout ouvert $U$ de $M$, si $f\vert_U$ est une fonction lisse, alors la distribution $u\vert_U$ est également une fonction lisse.

Cette condition de régularité est notamment satisfaite par les opérateurs elliptiques, tels que le Laplacien. Dans son célèbre théorème des sommes de carrés, \textcite{Hor67} a montré que si $\{X_1, \dots, X_n\}$ est une famille de champs de vecteurs satisfaisant la condition dite des crochets de Hörmander, alors l'opérateur différentiel $X_1^2 + \dots + X_n^2$ est hypoelliptique, bien qu'il ne soit généralement pas elliptique. Ce théorème a été le point de départ du développement de la théorie des opérateurs \emph{hypoelliptiques maximaux}, notamment à travers les travaux de \textcite{FS74, RS76, HN85}. 

La \emph{condition des crochets de Hörmander} porte sur l'existence d'un entier $N$ tel que tout vecteur tangent à $M$ en un point $x$ peut s'écrire comme une combinaison linéaire des évaluations en $x$ des champs $X_1, \dots, X_n$ et de leurs crochets de profondeur au plus $N$. On dit alors que la famille $\{X_1, \dots, X_n\}$ induit une \emph{structure sous-riemannienne} sur la variété $M$. On peut noter que sous ces conditions, tout opérateur différentiel peut s'écrire comme un polynôme non commutatif à coefficients dans $C^{\infty}(M)$ en $X_1, \dots, X_n$.  

\smallskip L'objet du travail présenté ici est l'étude de l'\emph{hypoellipticité maximale} des opérateurs différentiels obtenus comme polynômes en les éléments d'une famille de champs de vecteurs $\{X_1, \dots, X_n\}$ satisfaisant la condition des crochets de Hörmander. Cette notion plus fine étend la propriété de régularité qui caractérise les opérateurs différentiels hypoelliptiques en termes de régularité Sobolev, offrant ainsi un parallèle, en géométrie sous-riemannienne, aux opérateurs elliptiques.

\smallskip Dans leur article, \textcite{HN79c} ont proposé une conjecture caractérisant l’hypoellipticité maximale, généralisant le théorème principal de régularité des opérateurs elliptiques. Cette conjecture a été récemment confirmée grâce à des outils de géométrie non commutative. Un élément central de ce travail est une généralisation naturelle du groupoïde tangent de \textcite{ConnesNCG}, introduite en géométrie sous-riemannienne par \textcite{MoBlup,MoTanG}, dans lequel apparaissent tous les cônes tangents, ingrédients clés du travail de Helffer et Nourrigat. Ce travail s’appuie également sur la construction du calcul pseudodifférentiel à partir de l'action naturelle de $\mathbb{R}_+^*$ sur le groupoïde tangent \parencite{DS1,EY}.  

Dans leur article, \textcite{AMY} ont développé un calcul pseudodifférentiel dans ce contexte de géométrie sous-riemannienne, introduisant notamment la notion de symbole principal. Ils ont démontré que l’inversibilité de ce symbole équivaut à l’hypoellipticité maximale, validant ainsi la conjecture.

\medskip Il faudrait un livre pour expliquer en détail tous les prérequis et les résultats développés dans l'article de \textcite{AMY}, ainsi qu’une immense bibliographie pour évoquer tous les mathématiciens ayant contribué, d'une façon ou d'une autre, aux résultats présentés ici. Je me contenterai de donner les grandes lignes (d'une partie réduite) des travaux d'Androulidakis, Mohsen et Yuncken, et je prie par avance les nombreux auteurs que je n’ai pas mentionnés de bien vouloir m’excuser. 

L'article est organisé comme suit :

\begin{itemize}
    \item[\pt] La conjecture de Helffer et Nourrigat est présentée dans la première section.
    \item[\pt] La seconde section offre un aperçu de l'approche par les groupoïdes du calcul pseudodifférentiel classique.
    \item[\pt] La dernière section est consacrée à la description du groupoïde tangent dans le contexte de la géométrie sous-riemannienne et à la présentation des résultats clés permettant d’établir le théorème de régularité des opérateurs différentiels hypoelliptiques maximaux. L’article se termine par l'application de ces résultats à un exemple détaillé.
\end{itemize}

\smallskip Je me suis appuyée sur de nombreux documents fournis par Androulidakis, Mohsen et Yuncken pour rédiger ce texte, et je les remercie pour l’aide apportée. Cet article s’inspire fortement du mémoire d'habilitation de Mohsen, des discussions que nous avons eues et des nombreuses questions (petites ou grandes) auxquelles il a répondu. Enfin, je remercie Omar Mohsen et Georges Skandalis pour leur relecture attentive de ce document.

\subsection*{Avant-propos}

Les résultats obtenus par \textcite{AMY} sont très généraux et tolèrent en particulier l'adjonction de \enquote{poids} sur les champs de vecteurs $X_1, \dots, X_n$. Nous en présentons ici une version simplifiée en passant sous silence autant que possible les difficultés techniques. En particulier :

\begin{itemize}
    \item[\pt] Dans toute la suite, $M$ est une variété lisse compacte et $\cX(M)$ désigne le $C^{\infty}(M,\R)$-module des champs de vecteurs réels lisses sur $M$.
    \item[\pt] Nous considérerons des fonctions à valeurs complexes, mais, la plupart du temps, on pourrait également considérer des fonctions à valeurs dans un fibré (en particulier dans la définition des opérateurs pseudodifférentiels ou pour définir l'intégration sans ambiguïté).
    \item[\pt] Nous ne justifierons pas que les intégrales utilisées ont un sens. Suivant ses préférences, on pourra supposer :
        \begin{itemize}
            \item avoir choisi la mesure de Lebesgue si l'on travaille sur $\R^n$ ;
            \item avoir fait le choix d'une forme volume pour une variété orientable ;
            \item que les fonctions considérées prennent leurs valeurs dans un fibré de densité bien choisi dans le cas général.
        \end{itemize}
        Lorsque l'on n'utilise pas les densités, des termes correctifs liés aux changements de cartes doivent être pris en compte dans les formules.
    \item[\pt] Les espaces de Sobolev seront considérés avec des indices entiers, mais ils peuvent être étendus (ainsi que les résultats correspondants) aux indices réels via les techniques classiques de dualité et d'interpolation. Il en va de même pour l'ordre des opérateurs pseudodifférentiels.
\end{itemize}

\section{Opérateurs différentiels hypoelliptiques maximaux}

\smallskip

Dans cette section, nous donnons quelques rudiments de géométrie sous-riemannienne nécessaires à la suite. Nous rappelons ensuite le théorème de régularité des opérateurs différentiels elliptiques, puis présentons les opérateurs hypoelliptiques ainsi que la conjecture de Helffer et Nourrigat.

\subsection{Quelques mots sur les cônes tangents en géométrie sous-riemannienne}\label{conetang}

Dans toute la suite,  $\{X_1,\dots,X_n\}$ est une famille de champs de vecteurs sur $M$ satisfaisant la condition des crochets de Hörmander. On note $\cF=\cF^1$ le sous-module de $\cX(M)$ engendré par $X_1,\dots,X_n$, puis pour $i\geq1$, $\cF^{i+1}$ celui engendré par $\cF^i$ et $[\cF^1,\cF^i]$. La \emph{condition de Hörmander} dit précisément qu'il existe $N\in \N$ tel que $\cF^N=\cX(M)$. \\ Le \emph{sous-espace horizontal} en $x\in M$ est l'image $F_x:=\mbox{vect}(X_1(x),\dots,X_n(x))$ de $\cF$ dans $T_xM$.

\smallskip Ce cadre induit une métrique intrinsèque sur $M$, appelée métrique de \emph{Carnot-Carathéodory}, que nous décrivons ci-après. L'idée est de mesurer la distance entre deux points en suivant uniquement des chemins tangents aux sous-espaces horizontaux. On parle alors de \emph{structure sous-riemannienne} sur $M$. L'essentiel de ce qui est présenté dans cette partie provient de l'article de \parencite{Bel}.

\smallskip Pour $x \in M$ et $v \in F_x$, on définit 
$$\|v\|:= \inf \{\sqrt{\sum_{i=1}^n v_i^2} \ ; \ v=\sum_{i=1}^n v_iX_i(x)\}$$

Un chemin absolument continu $\gamma :[0,1] \rightarrow M$ est dit \emph{horizontal} si pour presque tout 
$t$, $\gamma'(t) \in F_{\gamma(t)}$. Sa \emph{longueur} est donnée par la formule
$$\ell(\gamma)=\int_0^1 \|\gamma'(t)\| dt$$

\medskip Le théorème de Chow-Rashevskii assure que deux points de $M$ sont toujours reliés par un chemin horizontal. On définit alors :
\begin{defi}\label{dcc}
La distance de \emph{Carnot-Carathéodory} entre les points $x$ et $y$ de $M$ est 

\centerline{$d_{CC}(x,y)=\inf \{ \ell(\gamma) \ ; \ \gamma \mbox{ chemin horizontal tel que } \gamma(0)=x \mbox{ et  } \gamma(1)=y\}$}
\end{defi}
Notons que cette distance engendre la topologie usuelle de $M$.

\smallskip 
Rappelons que la \emph{distance de Gromov-Hausdorff} mesure à quel point deux espace métriques compacts sont non isométriques \footnote{Pour deux espace métriques $A$ et $B$, $d_{HG}(A,B)=\inf \{d_H(f(A),g(B))\}$ où $(f,g,Z)$ parcourent l'ensemble des couples de plongements isométriques $f:A\rightarrow Z$, $g:B\rightarrow Z$ et $d_H$ est la distance de Hausdorff.}. \\ On dira qu'un espace métrique $(T,d)$ admet un \emph{cône tangent} en $x\in T$ si les boules $B(x,\varepsilon)$ munies de la distance $d/\varepsilon$, convergent au sens de Gromov-Hausdorff lorsque $\varepsilon$ tend vers $0$. 

\smallskip Par exemple, si $M$ est une variété riemannienne, son cône tangent en un point $x$ est son espace tangent $T_xM$ avec sa métrique euclidienne.

\smallskip On sait depuis la généralisation des travaux de \textcite{Mit} par \textcite{Bel} que, dans le cas d'une variété sous-riemannienne $M$ munie de la métrique de Carnot-Carathéodory $d_{CC}$, la limite en $0$ au sens de Gromov-Hausdorff pointée de $(M,d_{CC}/\varepsilon,x)$ est unique et qu'il s'agit d'un espace homogène d'un groupe de Lie nilpotent. 

\smallskip Précisément, notons $\frak{g}$ l'algèbre de Lie libre engendrée par $n$ éléments notés $\tilde{X_1},\cdots,\tilde{X_n}$ ayant pour relations que les itérations de crochets de longueur strictement plus grande que $N$ s'annulent. Il s'agit d'une algèbre de Lie nilpotente de dimension finie. Notons $G$ le groupe de Lie nilpotent simplement connexe qui l'intègre.

\smallskip \textcite{Bel} exhibe, pour tout $x$ de $M$, un sous-groupe simplement connexe $G_x$ de $G$ de codimension $\mbox{dim} (M)$ tel que :
$$\underset{\varepsilon \rightarrow 0}{\mbox{lim}} (M, d_{CC}/\varepsilon,x)=(G/G_x, d_{G/G_x},G_x)$$ où $d_{G/G_x}$ est une métrique de Carnot-Carathéodory sur l'espace homogène $G/G_x$.

\medskip Les travaux de \textcite{MoTanG} révèlent qu'il existe d'autres cônes tangents à considérer et nous verrons dans la section \ref{TangOmar} comment les obtenir tous.

\subsection{Opérateurs différentiels hypoelliptiques maximaux}

Une \emph{distribution} sur la variété compacte lisse $M$ est un élément du dual topologique de $C^{\infty}(M)$, c'est-à-dire une forme linéaire $u:C^{\infty}(M)\rightarrow \C$ continue pour la topologie naturelle : il existe $k\in \N$ tel que $u$ s'étend en une forme linéaire continue sur l'espace de Banach $C^k(M)$. Le plus petit tel $k$ est appelé l'\emph{ordre} de la distribution $u$. On note $C^{-\infty}(M)$ l'espace des distributions sur $M$.

\smallskip 
Le choix d'une $1$-densité $dx$ sur $M$ permet de plonger $C^{\infty}(M)$ dans $C^{-\infty}(M)$ : la fonction $f\in C^{\infty}(M)$ définit la distribution $\begin{aligned} g\in C^{\infty}(M)\mapsto \int_M f(x)g(x)dx \end{aligned}$.

\smallskip  Un opérateur différentiel $D$ sur $M$ s'étend  naturellement en un endomorphisme linéaire encore noté $D$ de l'espace des distributions $C^{-\infty}(M)$. 

\smallskip Un tel opérateur $D$ est dit \emph{hypoelliptique} si, pour tout ouvert $U$ de $M$ et toute fonction lisse $f$ en restriction à $U$, la solution $u$ de l'équation $Du=f$ est elle aussi lisse sur $U$.

\begin{exem} On se place sur $\R^2$.
\begin{enumerate} \item Le \emph{Laplacien} : $P=\partial_x^2+\partial_y^2$ est hypoelliptique.
\item L'opérateur de \emph{la chaleur} : $P=\partial_t-\partial_x^2$ est hypoelliptique.
\item L'opérateur \emph{d'Alembertien} associé à l'équation des ondes : $P=\partial_t^2-\partial_x^2$ n'est pas hypoelliptique.
\item Le \emph{Laplacien hypoelliptique} : $P=-\partial_x^2-x\partial_y$ est hypoelliptique.
\item Le sous-Laplacien de Grushin : $P_{\lambda}=-\partial_x^2-x^2\partial_y^2+\lambda \partial_y$ est hypoelliptique si et seulement si $\lambda \notin \{(2k+1)i, k\in \Z\}$. 
\end{enumerate}
\end{exem}

\subsubsection*{\Gt Rappel sur le théorème de régularité des opérateurs différentiels elliptiques}$\ $

\smallskip
Le Laplacien est un exemple d'opérateur différentiel elliptique. Les opérateurs différentiels elliptiques bénéficient de propriétés de régularités très fortes énoncées plus bas qui garantissent, en particulier qu'ils sont hypoelliptiques. 

\medskip Rappelons qu'un opérateur différentiel $D: C^{\infty}(M) \rightarrow C^{\infty}(M)$ d'ordre $k$ s'écrit localement $\sum_{\vert \alpha \vert \leq k} f_{\alpha} \partial^{\alpha}$ où  si $m$ désigne la dimension de la variété $M$, $\alpha=(\alpha_1,\cdots,\alpha_m) \in \N^m$ est un multi-indice de longueur $|\alpha |=\sum \alpha_i$ et $\partial^{\alpha}=\dfrac{\partial^{\alpha_1}}{(\partial x_1)^{\alpha_1}}\cdots\dfrac{\partial^{\alpha_m}}{(\partial x_1)^{\alpha_m}}$.

On définit localement le \emph{symbole principal} de $D$ par $$\sigma(D,x,\xi)=\sum_{\vert \alpha \vert = k} f_{\alpha}(x)i^k\xi_1^{\alpha_1}\cdots\xi_m^{\alpha_m}$$
On vérifie que cette formule est indépendante du choix de coordonnées et que le symbole principal défini une fonction lisse $\sigma(D):T^*M\rightarrow \C$ homogène de degré $k$.

\smallskip Les résultats suivant s'obtiennent à l'aide du calcul pseudodifférentiel classique \parencite{HorV3}. 
\begin{theo}\label{ell}
Soit $D$ un opérateur différentiel d'ordre $k$ sur la variété compacte et lisse $M$. Les assertions suivantes sont équivalentes:
\begin{enumerate}
\item L'opérateur $D$ est \emph{elliptique}, c'est-à-dire que pour tout $(x,\xi) \in T^*M\setminus M\times \{0\}$, $\sigma(D,x,\xi)\not=0$.
\item Pour tout (ou un) $s\in \N$ et toute distribution $u$ sur $M$, si $Du\in H^s(M)$ alors $u\in H^{s+k}(M)$.
\item Pour tout opérateur différentiel $P$ d'ordre inférieur ou égal à $k$, il existe $C>0$ tel que $$\|P(f)\|_{L^2(M)} \leq C\big(\|D(f)\|_{L^2(M)}+\|f\|_{L^2(M)}\big), \ \forall \ f\in C^{\infty}(M)$$
\item Pour tout (ou un) $s\in \N$, l'opérateur $D: H^{s+k}(M) \rightarrow H^{s}(M)$ est un opérateur de Fredholm \footnote{$M$ compact est nécessaire ici.}.
\end{enumerate}
\end{theo}

Puisque $C^{\infty}(M)=\underset{s\in \N}{\bigcap} H^s(M)$, on en déduit que si $D$ est elliptique, il est hypoelliptique.

\medskip Un ingrédient important dans la théorie des opérateurs elliptiques et qu'ils admettent des \emph{paramétrix}, c'est-à-dire des \enquote{inverses à régularisant près}. Un paramétrix pour l'opérateur $D$ est un opérateur $Q: C^{-\infty}(M)\rightarrow C^{-\infty}(M)$ tel que :
\begin{itemize}\item[\pt]$Q(C^{\infty}(M)) \subset C^{\infty}(M)$,
\item[\pt] les opérateurs $QD-Id$ et $DQ-Id$ sont \emph{régularisants}, c'est à dire : 

\centerline{$(QD-Id)(C^{-\infty}(M))\subset C^{\infty}(M)$ et $(DQ-Id)(C^{-\infty}(M))\subset C^{\infty}(M)$.}
\end{itemize}

\subsubsection*{\Gt Théorème de régularité des opérateurs différentiels hypoelliptiques maximaux et  conjecture de Helffer et Nourrigat} $\ $

\smallskip
Considérons à nouveau la famille de champs de vecteurs $\{X_1,\dots,X_n\}$ sur $M$ satisfaisant la condition des crochets de Hörmander. Cette dernière condition garantit que tout opérateur différentiel $D$ peut s'obtenir sous la forme $D=P(X_1,\cdots,X_n)$ où $P$ est un polynôme non commutatif à coefficients dans $C^{\infty}(M)$. \\ Nous renvoyons aux travaux de \textcite{HN85} pour les détails concernant cette partie.

\begin{defi} \label{SobetOp}\begin{enumerate} \item L'\emph{ordre de Hörmander} d'un opérateur différentiel $D$ est le plus petit degré de polynôme $P$ tel que $D=P(X_1,\cdots,X_n)$.  On notera $DO^k(\cF^{\bullet})$ l'ensemble des opérateurs différentiels d'ordre de Hörmander $k$.
\item On considère les espaces de Sobolev définis pour $s\in \N$ par :
$$\tilde{H}^s(M):=\{u\in L^2(M) \ : \ Du\in   L^2(M) \mbox{ pour tout } D \mbox{  d'ordre de Hörmander } \leq s\}$$

\end{enumerate}
\end{defi}

On a encore $C^{\infty}(M)=\underset{s\in \N}{\bigcap}\tilde{H}^s(M)$.

\begin{exem} Dans $\R^2$, avec $X_1=\partial_x$ et $X_2=x\partial_y$, l'ordre de Hörmander de $\partial_y$ est $2$ car $\partial_y=[X_1,X_2]$.
\end{exem}

\medskip Il s'agit maintenant de définir un symbole principal dans ce cadre. \\
On considère à nouveau $\frak{g}$, l'algèbre de Lie libre engendrée par $n$ éléments notés $\tilde{X_1},\cdots,\tilde{X_n}$ ayant pour relations que les itérations de crochets de longueur strictement plus grande que $N$ s'annulent. On note $G$ le groupe de Lie nilpotent simplement connexe qui l'intègre et  $\Xi : \frak{g} \rightarrow G$ l'application exponentielle, qui est un difféomorphisme ici.

\smallskip Soit $\pi : G\rightarrow U(H)$ une représentation unitaire irréductible de $G$.

On note $C^{\infty} (\pi)$ le sous-espace de $H$ des \emph{vecteurs lisses}, c'est-à-dire des vecteurs $v\in H$ tels que $g\in G \mapsto \pi (g)v \in H$ est lisse.

\smallskip Rappelons que la différentielle de $\pi$ : $$d\pi : \frak{g} \rightarrow \operatorname{End}\big(C^{\infty}(\pi)\big)$$ est donnée par la formule $d\pi(X)v=\underset{t\rightarrow 0}{\mbox{lim }} \dfrac{\pi(\Xi(tX))v-v}{t}$.

\medskip Considérons un opérateur différentiel $D\in DO^k(\cF^{\bullet})$ et un polynôme $P$ de degré $k$ tel que $D=P(X_1,\cdots,X_n)$. On définit :
$$\sigma(D,x,\pi)=P_{\max,x}(d\pi(\tilde{X_1}),\cdots,d\pi(\tilde{X_n}) )\in \operatorname{End}\big(C^{\infty}(\pi)\big)$$ 
où $P_{\max,x}$ est obtenu en ne gardant que la partie homogène de degré $k$ de $P$ et en évaluant les fonctions coefficients en $x$.

\medskip À priori, $\sigma(D,x,\pi)$ dépend du choix du polynôme $P$. Dans leur article, \textcite{HN79c} construisent, à l'aide de la méthode des orbites de Kirillov, un ensemble $\cT_x$ de représentations unitaires, appelé \emph{ensemble de Helffer-Nourrigat}, qui ne dépend que des champs de vecteurs $\{X_1,\dots,X_n\}$ et pour lequel l'application ci-dessus ne dépend pas du choix du polynôme de degré $k$ qui définit l'opérateur $D$. On donnera une description des $\cT_x$ en \ref{Tau} et leur lien avec les cônes tangents.

\medskip On peut énoncer le théorème de régularité des opérateurs différentiels hypoelliptiques maximaux \parencite{AMY} :

\begin{theo}\label{hypo} Soit une famille de champs de vecteurs $\{X_1,\dots,X_n\}$ sur la variété compacte et lisse $M$ satisfaisant la condition des crochets de Hörmander et $D$ un opérateur différentiel sur $M$ d'ordre de Hörmander $k$. 

 Les assertions suivantes sont équivalentes:
\begin{enumerate}
\item Pour tout $x\in M$ et $\pi\in \cT_x\setminus\{1_G\}$, $\sigma(D,x,\pi)$ est injective.
\item Pour tout (ou un) $s\in \N$ et toute distribution $u$ sur $M$, si $Du\in \tilde{H}^s(M)$ alors $u\in \tilde{H}^{s+k}(M)$.
\item Pour tout opérateur différentiel $P$ d'ordre de Hörmander inférieur ou égal à l'ordre de Hörmander de $D$, il existe $C>0$ tel que 
$$\|P(f)\|_{L^2(M)} \leq C\big(\|D(f)\|_{L^2(M)}+\|f\|_{L^2(M)}\big), \ \forall \ f\in C^{\infty}(M)$$
\item Pour tout (ou un) $s\in \N$, l'opérateur $D: \tilde{H}^{s+k}(M) \rightarrow \tilde{H}^{s}(M)$ est inversible à gauche modulo les opérateurs compacts \footnote{$M$ compact est nécessaire ici.}.
\end{enumerate}
\end{theo}

Un opérateur différentiel satisfaisant ces conditions est dit \emph{hypoelliptique maximal}. 
\smallskip Puisque $C^{\infty}(M)=\underset{s\in \N}{\bigcap} \tilde{H}^s(M)$, un tel opérateur est en particulier hypoelliptique.

\medskip Ce théorème était essentiellement conjecturé par \textcite{HN79c} jusqu'aux travaux de \textcite{AMY}. Précisément, on peut noter :
\begin{itemize} \item L'implication de (2) vers (1) a été montrée par \textcite{HN85} en toute généralité.
\item L'implication de (1) vers (2) a été montrée dans certains cas particuliers (rang 2 ou équirégulier) par \textcite{RS76}, \textcite{Rot79}, \textcite{HN85}.
\item \textcite{Str14} montre que si l'inégalité (3) est vraie pour la norme $L^2$ alors elle est vraie pour n'importe quelle norme $L^p$, $p\in ]1,+\infty[$.
\item Les résultats de Helffer et Nourrigat sont obtenus sans paramétrix. De nombreux auteurs ont construit des calculs pseudodifférentiels dans différentes généralités qui contiennent les paramétrix mais pour lesquels le symbole principal est difficile à calculer. Le livre de \textcite{HN85} fournit un panorama de ces travaux ainsi qu'une bibliographie complète. Concernant des travaux plus récents, citons : \textcite{BG88,  EMM91, EMM91, Chr92, Pon08, BFG12, FR14, FR16, EY, FF20, Cr24}.
\end{itemize}

\begin{rema} Si, pour tout $x\in M$, $F_x=vect(X_1(x),\cdots,X_n(x))=T_xM$ alors les théorèmes \ref{ell} et \ref{hypo} sont les mêmes.
\end{rema}

\section{Groupoïdes de Lie et calculs pseudodifférentiels classique}

Dans les années 80, deux contributions majeures presque simultanées ont imposé les groupoïdes dans le monde de la géométrie non commutative. D'une part, la construction par \textcite{Ren} de la $C^*$-algèbre d’un groupoïde localement compact. 
D'autre part, la construction par  \textcite{ConnesLNM,ConnesSurvey} de la $C^*$-algèbre d’un feuilletage, basée sur son groupoïde d’holonomie ainsi que d'un calcul pseudodifférentiel longitudinale \enquote{porté} par ce même groupoïde.
Suivra la construction du \emph{groupoïde tangent} de \textcite{ConnesNCG} qui permet 
de construire l’indice analytique des opérateurs (pseudo)différentiels sans utiliser le calcul pseudodifférentiel et de fournir une preuve alternativee du théorème d’indice d’Atiyah-Singer.

Le calcul pseudodifférentiel sur les groupoïdes, la construction de l’indice analytique et le groupoïde tangent de Connes ont inspiré de nombreux travaux où ces idées ont été généralisées à des contextes géométriques variés, nous renvoyons à la bibliographie très complète qui se trouve dans l'article de synthèse \parencite{DSreview}.

Plus récemment, nous avons montré avec Skandalis comment obtenir directement les opérateurs pseudodifférentiels en utilisant l'action canonique de $\R_+^*$ sur le groupoïde tangent \parencite{DS1}. Une construction alternative reposant sur la même idée a été ensuite développée par \textcite{EY}.

\smallskip  Dans ce préambule, nous rappelons les définitions de groupoïdes de Lie, de leurs $C^*$-algèbres et des opérateurs pseudodifférentiels associés.

\subsection{Groupoïdes}

\begin{defi}
Un \emph{groupoïde} est une petite catégorie dont toutes les flèches sont inversibles. Autrement dit, un \emph{groupoïde} est donné par :
\begin{itemize}
\item[\pt] Un ensemble d'\emph{unités} (d'objets) $G^{(0)}$, un ensemble de  \emph{flèches} (morphismes) $G^{(1)}=G$.
\item[\pt] Une application unité injective $u:G^{(0)} \rightarrow G$.\\ On identifie $G^{(0)}$ à son image dans $G$.
\item[\pt] Des applications \emph{source} et  \emph{but} $r,s:G\rightarrow G^{(0)}$ telles que $s\circ u=r\circ u=Id$.
\item[\pt] Un \emph{produit} partiel $G^{(2)}:=\{(\gamma,\eta)\in G \ \vert \ s(\gamma)=r(\eta)\} \rightarrow G,\ (\gamma,\eta)\mapsto \gamma\cdot \eta$ qui vérifie :
\begin{itemize} \item $s(\gamma\cdot \eta)=s(\eta)$ et $r(\gamma\cdot \eta)=r(\gamma)$,
\item le produit est associatif : $(\gamma_1\cdot \gamma_2)\cdot \gamma_3=\gamma_1\cdot (\gamma_2\cdot \gamma_3)$ dès que $s(\gamma_1)=r(\gamma_2)$ et $s(\gamma_2)=r(\gamma_3)$,
\item compatible avec les unités : $\gamma\cdot s(\gamma)=r(\gamma)\cdot \gamma=\gamma$.
\end{itemize} 
\item[\pt] Une application \emph{inverse} : $G\rightarrow G, \gamma \mapsto \gamma^{-1}$ qui est involutive et satisfait :
\begin{itemize} 
\item $s(\gamma^{-1})=r(\gamma)$, 
\item $\gamma \cdot \gamma^{-1}=r(\gamma)$. 
\end{itemize} 
\end{itemize} 

\smallskip On note alors $G \rightrightarrows G^{(0)}$.
\end{defi}

\begin{rema} Le groupoïde induit une relation d'équivalence sur l'espace des unités $G^{(0)}$ :

\centerline{pour $x,\ y$ dans $G^{(0)}$ : $x\sim y \iff \exists \gamma \in G \ :\ s(\gamma)=x \mbox{ et } r(\gamma)=y$.}

\medskip Pour $x\in G^{(0)}$, l'\emph{orbite} de $G$ passant par $x$ est $\cO_x=r(s^{-1}(x))=s(r^{-1}(x))$. On note $G^{(0)}/G$ l'\emph{espace des orbites}.

\smallskip Le \emph{groupe d'isotropie} en $x$ : $G_x^x:=s^{-1}(x)\cap r^{-1}(x)$ est un groupe d'unité $x$. \\
Le groupe $G_x^x$ agit sur $G_x:=s^{-1}(x)$ par multiplication à droite et sur $G^x:=r^{-1}(x)$ par multiplication à gauche :
$$G_x/G_x^x \simeq \cO_x \simeq G^x/G_x^x$$
\end{rema}

\medskip On parle de \emph{groupoïde localement compact} lorsque $G^{(0)}$ et $G$ sont des espaces topologiques tels que :
\begin{itemize} \item[\pt] $G$ et $G^{(0)}$ sont des espaces localement compacts à base dénombrable,
\item[\pt] toutes les applications de structure (unité, source but, produit, inverse) sont continues et la source et le but sont ouvertes.
\end{itemize}

\medskip Un groupoïde localement compact $G\rightrightarrows G^{(0)}$ est un \emph{groupoïde de Lie} lorsque $G$ et $G^{(0)}$ sont des variétés lisses, $G^{(0)}$, $s$ et $r$ sont des submersions surjectives (donc $G^{(2)}$ est une variété), le produit et l'inverse sont lisses et $u$ est un plongement. 

\begin{rema} La définition usuelle de groupoïde de Lie autorise la variété $G$ a être une variété non séparée. Cette situation est fréquente lorsque l'on s'intéresse aux groupoïdes de feuilletages.  Les constructions de $C^*$-algèbres et d'opérateurs pseudodifférentiels qui sont présentées dans la suite s'adaptent au cas non séparé. Cependant nous ne considérerons pas ce cas ici, tous les groupoïdes rencontrés étant séparés.
\end{rema}

\begin{exem}
\begin{enumerate}
\item $M$ un ensemble est un groupoïde, $M\rightrightarrows M$ avec $s=r=Id$, $M^{(2)}=M$. Une variété est un groupoïde de Lie. 
\item $G$ un groupe d'unité $e$ est un groupoïde : $G\rightrightarrows \{e\}$, $G^{(2)}=G\times G$ et le produit est le produit usuel. C'est un groupoïde de Lie lorsque $G$ est un groupe de Lie.
\item Un fibré en groupes $\pi:E\rightarrow M$ est un groupoïde : $E\rightrightarrows M$, $s=r=\pi$, $E^{(2)}=\sqcup_{x\in M}E_x\times E_x \rightarrow E$ avec le produit standard. En particulier, un fibré vectoriel lisse est un groupoïde de Lie.
\item $M$ un ensemble, le \emph{groupoïde des couples} est $M\times M \rightrightarrows M$ avec 
$$s(x,y)=y,\ r(x,y)=x, \ (x,y)\cdot(y,z)=(x,z),\ u(x)=(x,x),\ (x,y)^{-1}=(y,x)$$
L'espace des orbites est réduit à un point : $\cO_x=M$ pour tout $x\in M$.\\ C'est un groupoïde de Lie lorsque $M$ est une variété.

\smallskip Plus généralement, le graphe d'une relation d'équivalence $\cR$ sur $M$  :
$$G_{\cR}=\{(x,y)\in M\times M \ ; \ x\cR y\}\rightrightarrows M$$ avec la même structure que précédemment est un groupoïde. \\ Dans ce cas, l'espace des orbites est $M/\cR$.
\item \emph{Action de groupe}

Soit $H$ un groupe de neutre $e$ agissant sur un ensemble $M$ à gauche, le groupoïde de l'action est : $H\ltimes M \rightrightarrows M$ où :

$H\ltimes M =H\times M$, $s(h,x)=x$, $r(h,x)=h.x$, $u(x)=(x,e)$, $(g,hx)\cdot (h,x)=(gh,x)$ et $(h,x)^{-1}=(h^{-1},hx)$.

Les orbites, l'espace des orbites et l'isotropie de $H\ltimes M$ coïncident avec ceux de l'action de $H$ sur $M$.
Une fois encore, si $H$ est un groupe de Lie, $M$ une variété lisse et que l'action est lisse, le groupoïde $H\ltimes M \rightrightarrows M$ est de Lie.

\item Soit $M$ une variété.
$$\Pi(M)=\{[\gamma] \ ; \ \mbox{ classe d'homotopie à extrémités fixes de } \gamma:[0,1]\rightarrow M \mbox{ continu } \} \rightrightarrows M$$
\begin{itemize}
\item $s([\gamma])=\gamma(0)$, $r([\gamma])=\gamma(1)$,
\item $u(x)=[C_x]$ où $C_x$ est le chemin constant en $x$,
\item $[\gamma]\cdot [\eta]=[\gamma \cdot \eta]$ lorsque $\gamma(0)=\eta(1)$ où $\gamma \cdot \eta=\left\{ \begin{array}{l} \eta(2t) \mbox{ pour } t\in [0,\dfrac{1}{2}] \\ \gamma(2t-1) \mbox{ pour } t\in [\dfrac{1}{2},1] \end{array} \right.$
\item $[\gamma]^{-1}=[\gamma^{-1}]$ où $\gamma^{-1}(t)=\gamma(1-t)$.
\end{itemize}
\smallskip On munit naturellement  $\Pi(M)$ d'une structure de variété de dimension $2\mbox{dim}M$ pour laquelle $\Pi(M)$ est un groupoïde de Lie, $\Pi(M)$ est un revêtement de $M\times M$.

\smallskip Pour $x\in M$, $\cO_x$ est la composante connexe de $x$, 
$\Pi(M)_x^x=\pi_1(M,x)$ est le groupe d'homotopie en $x$ et $\Pi(M)_x$ un revêtement universel.
\end{enumerate}
\end{exem}

Le pendant infinitésimal des groupoïdes de Lie sont les \emph{algébroïdes de Lie}.
\begin{defi} Un \emph{algébroïde de Lie} $(\cA,\natural,[\ ,\ ])$ sur une variété $M$ est un fibré vectoriel $\cA\rightarrow M$ tel que l'espace de ses sections lisses $(\Gamma(\cA),[\ ,\ ])$ est une algèbre de Lie  et  $\natural : \cA \rightarrow TM$ est un morphisme de fibré appelé \emph{ancre}, dont l'application induite sur les sections des fibrés $\natural : \Gamma(\cA) \rightarrow \cX(M)$ satisfait :
\begin{enumerate}
\item $\natural( [X,Y] )=[\natural(X),\natural(Y)]$,
\item $[X,fY]=f[X,Y]+(\natural(X)f)Y$,
\end{enumerate}
lorsque $X$ et $Y$ sont deux sections lisses de $\cA$ et $f$ est une fonction lisse sur $M$.
\end{defi}

\smallskip 

Si $G  {\rightrightarrows}G^{(0)}$  est un groupoïde de Lie, notons qu'il peut être considéré comme la famille paramétrée par $G^{(0)}$ des variétés $G_x:=\{\gamma \in G \ ; \ s(\gamma)=x\},\ x\in G^{(0)}$.

La \emph{translation à droite} donnée pour tout $\gamma \in G$ avec $y=r(\gamma)$ et $x=s(\gamma)$, par 
$R_{\gamma} : G_y \rightarrow G_x,\ \eta \mapsto \eta\cdot \gamma$ agit sur la famille $(G_x)_{x\in G^{(0)}}$.

\smallskip On construit l'algébroïde de Lie du groupoïde $G\rightrightarrows G^{(0)}$ de la façon suivante (voir \cite{Mack}). 
\begin{itemize}
\item[\pt] Le fibré vectoriel ${\gA}  G \rightarrow G^{(0)}$ est la restriction aux unités, $G^{(0)}$, du noyau de la différentielle $Ts$ de l'application source $s$. Ainsi ${\gA}  G=\cup_{x\in G^{(0)}} T_xG_x$ est la réunion des espaces tangents en les unités aux fibres de $s$.
\item[\pt] L'ancre est la restriction de la différentielle $Tr$ de l'application but $r$.
\item[\pt] Si $Z$ est une section de $\cA$, on lui associe un \emph{champ invariant à droite}, $\tilde{Z}$ défini par 
$$\widetilde{Z}(\eta)=TR_{\eta}(Z(r(\eta))$$
où $R_{\eta} : G_{r(\eta)} \rightarrow G_{s(\eta)}$ est la multiplication à droite par $\eta$.

Le crochet de Lie est alors donné par :
$$\begin{array}{cccc} [\ ,\ ]: & \Gamma({\gA}  G)\times \Gamma({\gA} 
  G) & \longrightarrow & \Gamma({\gA}  G) \\
 & (Z_1,Z_2) & \mapsto & [\widetilde{Z_1},\widetilde{Z_2}]_{G^{(0)}}
\end{array}$$
où $[\widetilde{Z_1},\widetilde{Z_2}]_{G^{(0)}}$  désigne la restriction du crochet usuel $[\widetilde{Z_1},\widetilde{Z_2}]$ aux unités ${G^{(0)}}$.
\end{itemize}
Par construction l'image de $\Gamma({\gA}  G)$ par l'ancre est la réunion des espaces tangents aux orbites de $G$.

\begin{exem}
\begin{enumerate}
\item L'algébroïde de Lie d'un groupe de Lie est l'algèbre de Lie du groupe.
\item L'algébroïde de Lie du groupoïde des couples $M\times M \rightrightarrows M$ sur une variété $M$ est $TM$ avec l'identité pour ancre.
\item Lorsque $\cF$ est un feuilletage régulier sur la variété $M$, $T\cF$ avec l'identité pour ancre est un algébroïde de Lie. Le \emph{groupoïde d'holonomie} de $\cF$ est le \enquote{plus petit} groupoïde de Lie dont $T\cF$ est  l'algébroïde de Lie. 
\end{enumerate}
\end{exem}

\subsubsection{$C^*$-algèbre}

Dans la suite $G\rightrightarrows G^{(0)}$ est un groupoïde de Lie (séparé).

\paragraph*{\Gt Systèmes de Haar}
On considère un analogue de la mesure de Haar des groupes localement compacts.

\begin{defi} Un \emph{système de Haar} sur le groupoïde de Lie $G\rightrightarrows G^{(0)}$ est une famille $\mu=\{\mu_x\ ; \ x\in G^{(0)} \}$ de mesures de Lebesgue $\mu_x$ sur $G_x$ satisfaisant les conditions suivantes :
\begin{description} \item[Support ] Pour tout $x\in G^{(0)}$, le support de $\mu_x$ est $G_x$.
\item[Invariance ] Pour tout $\gamma \in G$, $y=r(\gamma)$ et $x=s(\gamma)$, l'opérateur de translation $R_{\gamma} : G_y \rightarrow G_x$ préserve la mesure. C'est-à-dire, pour tout $f\in C^{\infty}_c(G)$ : 
$$\int_{G_x}f(\eta) d\mu_x(\eta)=\int_{G_y}f(\eta\cdot \gamma) d\mu_y(\eta)$$
\item[Différentiabilité  ] Pour tout $f\in C^{\infty}_c(G)$, l'application suivante est lisse : 
$$G^{(0)}\rightarrow \C \ ; \ x\mapsto \int_{G_x}f(\gamma) d\mu_x(\gamma)$$
\end{description}
\end{defi}

\paragraph*{\Gt Algèbre de convolution}
\bigskip Munis d'un système de Haar $\mu$ sur le groupoïde $G$, on équipe $C^{\infty}_c(G)$, l'espace vectoriel des fonctions lisses à support compact sur $G$, d'une structure d'algèbre de convolution donnée par :

\begin{description}
\item[Involution] Pour $f\in C^{\infty}_c(G)$, $\gamma \in G$, $f^*(\gamma)=\overline{f(\gamma^{-1})}$.
\item[Produit de convolution] Pour $f\in C_c^{\infty}(G)$, $\gamma \in G$ et $x=s(\gamma)$ 
$$f\ast g(\gamma)=\int_{\eta \in G_x}f(\gamma\cdot \eta^{-1})g(\eta) d\mu_x(\eta)$$
\end{description}

\paragraph*{\Gt Normes et représentations}

On définit la norme suivante sur $C^{\infty}_c(G)$ : 
$$\|f\|_I=\mbox{sup}_{x\in G^{(0)}}\{\max \big(\int_{G_x} |f(\gamma) | d\mu_x(\gamma), 
\int_{G_x} |f(\gamma^{-1}) | d\mu_x(\gamma) \big)\}$$

\begin{defi} 
Une représentation \emph{non dégénérée} de $C^{\infty}_c(G)$ est un morphisme d'algèbres involutives $\pi :C^{\infty}_c(G) \rightarrow \cL(H)$ où $H$ est un espace de Hilbert tel que l'ensemble $\{\pi(f)\cdot \xi \ ; \ f\in C_c(G),\ \xi \in H\}$ est dense dans $H$.

\smallskip On dit que $\pi$ est \emph{bornée} lorsque pour tout $f\in C^{\infty}_c(G)$, $\|\pi(f)\| \leq \|f\|_I$.
\end{defi}
 On définit la \emph{norme maximale} de $f\in C^{\infty}_c(G)$, par $$\|f\|_{\max}= \underset{\pi \ bornée}{\mbox{sup}} \| \pi(f) \|_{\cL(\cH)}$$

On peut définir une autre norme en considérant la \emph{famille des représentations régulières} :

pour $x \in G^{(0)}$, $\pi_x:C^{\infty}_c(G) \rightarrow \cL(L^2(G_x,\mu_x))$ est donnée par 
$$\pi_x(f)(\xi)(\gamma)=\int_{\eta \in G_x} f(\gamma \cdot \eta^{-1}) \xi(\eta)d\mu_x(\eta)$$

Pour tout $x\in G^{(0)} $ et $f\in C^{\infty}_c(G)$,  $\pi_x(f)$ est borné (i.e. dans $\cL(L^2(G_x,\mu_x))$ et $\pi_x$ est une représentation non dégénérée et bornée de $C^{\infty}_c(G)$.

On définit la \emph{norme minimale} de $f\in C^{\infty}_c(G)$, par $$\|f\|_{\min} = \underset{x\in G^{(0)} }{\mbox{sup}} \| \pi_x (f) \|_{\cL((L^2(G_x,\mu_x))}$$

\smallskip

La \emph{$C^*$-algèbre réduite}  de $G$ est la complétion de $C^{\infty}_c(G)$ pour la norme minimale : $C_r^*(G)=\overline{C_c^{\infty}(G)}^{\min}$. La \emph{$C^*$-algèbre maximale}  de $G$ est la complétion de $C^{\infty}_c(G)$ pour la norme maximale $C^*(G)=\overline{C_c^{\infty}(G)}^{\max}$. 

\smallskip Les groupoïdes que nous considérerons dans la suite sont \emph{moyennables}, dans ce cas les $C^*$-algèbres réduites et maximales sont égales.

\begin{exem}\begin{enumerate}
\item Le groupoïde des couples $M\times M \rightrightarrows M$ est moyennable et sa $C^*$-algèbre est l'algèbre $\cK(L^2(M))$ des opérateurs compacts sur $L^2(M)$.
\item Un fibré vectoriel lisse $E\rightarrow M$  est un groupoïde moyennable et sa $C^*$-algèbre est isomorphe via la transformée de Fourier à l'algèbre $C_0(E^*)$ des fonctions continues nulles à l'infini sur $E^*$.
\item Dans le cas d'un groupe de Lie, on retrouve les notions usuelles de $C^*$-algèbres (maximale ou réduite) de groupe.
\end{enumerate}
\end{exem}

\begin{rema} Il est possible de se soustraire au choix d'un système de Haar et d'obtenir une définition canonique de la convolution et de l'involution en remplaçant l'espace $C^{\infty}_c(G)$ par l'espace des sections à support compact d'un fibré (en droites) de demi-densités sur le groupoïde $G$. Choisir un système de Haar, c'est choisir une trivialisation de ce fibré. Nous conviendrons dans la suite, sans plus de détail, que l'intégration le long des fibres (de la source ou du but) d'un groupoïde de Lie a un sens.
\end{rema}

\paragraph*{\Gt Et si $G$ n'est pas tout à fait un groupoïde de Lie ?} $\ $

\smallskip Notons que la régularité essentielle dans la construction présentée plus haut est la différentiabilité des fibres $G_x=s^{-1}(x)$ du groupoïde $G$, lorsque $x\in G^{(0)}$. Ainsi, on peut appliquer cette construction aux \emph{familles continues de groupoïdes} (ou groupoïdes $C^{\infty,0}$) introduites par \textcite{Pater}. Un groupoïde localement compact $G\rightrightarrows G^{(0)}$ est $C^{\infty,0}$ lorsque $G_x$ est lisse pour tout $x\in G^{(0)}$ et la structure lisse varie continument avec $x$. En outre pour tout $\eta \in G$, la multiplication à droite $R_{\eta} : G_{r(\eta)} \rightarrow G_{s(\eta)}$ est un difféomorphisme et ces difféomorphismes varient continument avec $\eta$. 

Pour définir proprement de telles structure, on \enquote{remplace} les fonctions lisses par des fonctions de type $C^{\infty,0}$. 
 
 \begin{defi} Si $M$ et $N$ sont des espaces topologiques, $g:\R^n\times M \rightarrow \R^m$ et $h:M\rightarrow N$ deux applications continues, on dira que l'application $f=(g,h) : \R^n\times M \rightarrow \R^m \times N$ est $C^{\infty,0}$ lorsque pour tout $x\in M$, $g(\cdot,x)$ est lisse et toutes les différentielles partielles de $g$ en les variables de $\R^n$ sont continues.
 \end{defi}

\subsection{Groupoïde tangent et calcul pseudodifférentiel classique}

Historiquement, les opérateurs pseudodifférentiels sur un groupoïde ont été définis comme des familles d'opérateurs pseudodifférentiels (classiques) sur les fibres $G_x$ du groupoïde $G$ invariants par $G$. Ils ont été introduits par \textcite{ConnesLNM} dans le cadre des feuilletages réguliers, puis étendus à tout groupoïde de Lie indépendamment par  \textcite{MonthPie} et par \textcite{NWX}. 

La façon la plus naturelle de définir les opérateurs pseudodifférentiels sur un groupoïde est de les voir comme des distributions sur $G$ conormales à $G^{(0)}$ au sens de \textcite{HorV3}. La présentation que nous en faisons ici est principalement issue de notre article \parencite{DS1} dont on peut trouver une exposition détaillée dans la thèse de \textcite{Mahsa}. Pour une étude approfondie des distributions conormales dans le cadre des groupoïdes, nous renvoyons aux travaux de \textcite{LMV}.

\subsubsection{Déformation au cône normal, zooming action et éclatement}$\ $

Soit $\cM$ une variété lisse et $V$ une sous-variété. Notons $\cN_V^\cM$ le fibré normal de l'inclusion de $V$ dans $\cM$.

\paragraph*{\Gt Déformation au cône normal}

La \emph{déformation au cône normal} de $V$ dans $\cM$ est :
$$\operatorname{DNC}(\cM,V)=\cM\times \R^* \cup \cN_V^\cM\times \{0\}$$
On équipe $\operatorname{DNC}(\cM,V)$ d'une structure de variété lisse en demandant que :
\begin{itemize}
\item[\pt] la projection naturelle de $\operatorname{DNC}(\cM,V)$ sur $\cM\times \R$  est lisse,
\item[\pt] si $f:\cM\rightarrow \R$ est une fonction lisse sur $\cM$ qui s'annule sur $V$ alors la fonction $f^{dnc}:\operatorname{DNC}(\cM,V) \rightarrow \R$ donnée par :
$$\left\{ \begin{array}{l} f^{dnc}(x,t)=\dfrac{f(x)}{t} \mbox{ pour } (x,t)\in \cM\times \R^*, \\ \ f^{dnc}(x,X,0)=df(x,X) \mbox{ pour } (x,X,0)\in \cN_V^\cM\times \{0\}\end{array} \right.$$ est lisse.
\end{itemize}

\smallskip On peut également décrire cette structure de variété en choisissant une métrique riemannienne sur $\cM$ qui permet d'identifier $\cN_V^{\cM}$ à un voisinage tubulaire $\cU$ de $V$ dans $\cM$ via un difféomorphisme $\theta : \cN_V^{\cM} \rightarrow \cU$. On recolle alors les structures de variétés usuelles sur $\cM\times \R^*$ et sur $\cN_V^{\cM}\times \{0\}$ à l'aide de la carte :
$$\begin{array}{cccc} \Theta : & \cN_V^{\cM}\times \R & \rightarrow & \Theta(\operatorname{DNC}(\cM,V)) \\
& (x,X,t) & \mapsto & \left\{ \begin{array}{ll} (\theta(x,tX),t) & \mbox{ si } t\not=0 \\ (x,X,t) & \mbox{ si } t=0 \end{array} \right. \end{array}$$

\begin{rema} La construction $\operatorname{DNC}$ est fonctorielle.
\end{rema}

\begin{defi}
L'action naturelle de $\R^*$ sur $\cM\times \R^*$ s'étend en une action naturelle sur $\operatorname{DNC}(\cM,V)$ appelée  \emph{zooming action}:
$$  \begin{array}{ccc} \operatorname{DNC}(\cM,V)\times \R^* & \longrightarrow & \operatorname{DNC}(\cM,V) \\ (z,t,\lambda) & \mapsto & (z, \lambda t) \  \mbox{pour} \ t\not=0  
\\    (x,X,0,\lambda) & \mapsto & (x,  \frac{1}{\lambda} X,0) \  \mbox{pour} \ t=0  \end{array} 
$$
\end{defi}

\paragraph*{\Gt Eclatement} La zooming action est libre et propre en dehors de $V\times \R$ vu dans $\operatorname{DNC}(\cM,V)$ par le plongement naturel. On définit l'\emph{éclatement} de $V$ dans $\cM$ par :
$$\operatorname{Blup}(\cM,V)=\big( \operatorname{DNC}(\cM,V)\setminus V\times \R \big)/\R^*=\cM\setminus V \cup \mathbb{P}(\cN_V^\cM)\ $$
et l'\emph{éclatement sphérique} de $V$ dans $\cM$ par :
$$\operatorname{SBlup}(\cM,V)=\big( \operatorname{DNC}(\cM,V)\setminus V\times \R \big)/\R_+^*=(\cM\setminus V)_+\cup (\cM\setminus V)_-\cup \mathbb{S}(\cN_V^\cM)\ $$

Remarquons que l'application :
$$\begin{array}{ccc} \operatorname{DNC}(\cM,V)  & \longrightarrow  & 
\operatorname{DNC}(\cM\times \R,V\times \{0\})  \\ =\cM\times \R^* \cup \cN_V^\cM \times\{0\} & & =\cM\times \R \times \R^*\cup \cN_V^\cM \times\{0\}\\
(z,t) & \mapsto & (z,t,t) \ \mbox{pour} \ t\not=0  
\\    (x,X,0) & \mapsto & (x,  (X,1),0) \ \mbox{pour} \ t=0  \end{array} 
$$
a son image dans $\operatorname{DNC}(\cM\times \R,V\times \{0\})\setminus V\times \{0\}\times \R$ et qu'elle induit un plongement de $\operatorname{DNC}(\cM,V)$ dans $\operatorname{SBlup}(\cM\times \R,V\times \{0\})$.

\subsubsection{Distributions conormales} $\ $

On notera $\operatorname{DNC}(\cM,V)_+$ la sous variété à bord $\cM\times \R_+^* \cup \cN_V^\cM\times \{0\}$ de 
$\operatorname{DNC}(\cM,V)$.

\medskip L'espace $\cS_c(M,V)$ des \emph{fonctions de Schwartz} sur $\operatorname{DNC}(\cM,V)_+$ est l'espace des fonctions $f$ à support compact qui se prolongent par $0$ en une fonction lisse sur $\operatorname{SBlup}(\cM\times \R,V\times \{0\})$.

\medskip Le sous-espace $\cJ (\cM,V)\subset \cS_c(\cM,V)$ est l'espace constitué des fonctions $f=(f_t)_{t\in \R_+}$ telles que la famille $(f_t)_{t\in \R_+^*}$ s'annule à tout les ordres en $0$ comme distribution sur $\cM$, ce qui signifie que pour tout $g\in C^{\infty}(\cM)$, la fonction $\begin{aligned} t\mapsto \langle f_t,g \rangle=\int_\cM f_t(x)g(x)dx \end{aligned}$ s'étend en une fonction lisse sur $\R$ qui est nulle sur $\R_-$. 

\begin{exem} L'exemple fondamental à considérer est celui d'un fibré vectoriel réel lisse $q: E \rightarrow V$ muni d'une métrique. On identifie $V$  à la section nulle de $E$.
On a le diagramme commutatif suivant :

\[
  \begin{CD}
    E\times \R @>\Theta>\phantom{\text{blablabla}}> \operatorname{DNC}(E,V) \\
    @V{\varphi}VV @VViV\\
    \mathbb{S}_{E\times \R}\times \R @>>\tilde{\Theta}> \operatorname{SBlup}(E\times \R,V\times \{0\})
  \end{CD}
\]

Où $i$ est le plongement mentionné plus haut et $\Theta$ est le difféomorphisme  $$\begin{array}{ccc} E\times \R & \overset{\Theta}{\longrightarrow} &  \operatorname{DNC}(E,V) \\
(x,X,t) & \mapsto & \left\{ \begin{array} {ll} (x,tX,t) \mbox{ lorsque } t\not=0 \\ (x,X,0) \mbox{ lorsque } t=0 \end{array} \right. \end{array}$$
Le difféomorphisme $\tilde{\Theta}$ est  construit de façon analogue à $\Theta$ (pour le fibré $E\times \R\rightarrow V\times \{0\}$) et l'application $\varphi$ est  le plongement donné par $\varphi(x,X,t)=(x,\frac{1}{\sqrt{\|X\|^2+1}}(X,1),\sqrt{\|X\|^2+1}\ t)$.

\medskip Notons $S_c(E\times \R)$ l'espace des fonctions de Schwartz sur le fibré $E\times \R \rightarrow V\times \R$ à support compact dans la direction $V\times \R$ et $\cJ_{E\times \R}^+$ les fonctions  $f$ de $S_c(E\times \R)$ telles que pour tout $g\in C^{\infty}(E)$, la fonction $\begin{aligned} t\in \R^*_+\mapsto \langle f_t,g \rangle=\int_E f_t(x)g(x)dx \end{aligned}$ s'étend en une fonction lisse sur $\R$ qui s'annule sur $\R_-$.
\smallskip Cette dernière condition est équivalente 
à ce que pour tout $g\in C^{\infty}(E)$ la fonction $t\mapsto  \langle f_t,g \rangle$ s'annule à tous les ordres lorsque $t$ tends vers $0$ où encore que la transformée de Fourier $\hat f$ s'annule à tous les ordres sur la section nulle $V\times \{0\} \subset E^*\times \R$.

\smallskip Une fonction $f\in C^{\infty} ( \operatorname{DNC}(E,V))$ est dans $\cJ(E,V)$ si et seulement si $f\circ \Theta$ est dans $\cJ_{E\times \R}^+$.
\end{exem}

Notons pour tout $m\in \Z$, $\cI^m_c(\cM,V)$ l'espace des distributions conormales de \textcite{Hormandera} d'ordre $m$. On montre dans \parencite{DS1} le théorème suivant.

\begin{theo}\label{theodistrib} Pour tout $m \in \Z$, $\cI^m_c(\cM,V)$ est l'ensemble des distributions de la forme $\begin{aligned} \int_{\R^+} t^{-\dim M} k_t\dfrac{dt}{t^{1+m}} \end{aligned}$ où $k=(k_t)_{t\in \R_+}$ parcourt $\cJ (\cM,V)$ \footnote{Le facteur $ t^{-\dim M}$ n'apparaît pas si l'on considère des densités.}.
\end{theo}

\subsubsection{Groupoïde tangent}\label{pdoclas}

Soit $M$ une variété compacte et lisse, le \emph{groupoïde tangent} de $M$ \parencite{ConnesNCG} est la déformation au cône normal $\operatorname{DNC}(M\times M,\Delta M)_+$ où $M\times M \rightrightarrows M$ est le groupoïde des couples et $\Delta M=\{(x,x)\ ; x\in M\}$ l'espace de ses unités. La construction $\operatorname{DNC}$ étant fonctorielle, on obtient un groupoïde de Lie :
$$\T=M\times M \times \R^*_+ \sqcup TM\times \{0\} \rightrightarrows M\times M\times [0,+\infty)$$
où
\begin{itemize}
\item[\pt] $\T$ est la réunion des groupoïdes des couples $M\times M \times \{t\}\rightrightarrows M \times \{t\}$ pour $t\in \R^*_+$ et du fibré vectoriel (vu comme groupoïde) $TM\times \{0\}\rightrightarrows M\times \{0\}$.
\item[\pt] La structure de Lie est obtenue en recollant les structures de variété  usuelles sur $M\times M \times \R^*_+$ et sur $TM\times \{0\}$ comme précédemment. 

La topologie obtenue est telle qu'une suite $(y_n,x_n,t_n)_{n\in \N}$ de $M\times M \times \R_+^*$ converge vers $(x,X,0)\in TM$ lorsque la suite $(t_n)_{n\in \N}$ tend vers $0$, les suites $(x_n)_{n\in \N}$ et $(y_n)_{n\in \N}$ tendent vers $x$ et dans une carte la suite $(\dfrac{y_n-x_n}{t_n})_{n\in \N}$ tend vers $X$.\end{itemize}

\paragraph*{Calcul pseudodifférentiel} 

La structure de groupoïde permet de montrer que $S_c(\T):=S_c(M\times M, \Delta M)\subset C^{\infty}(\T)$ est une algèbre de convolution dont $\cJ:=\cJ(M\times M, \Delta M)$ est un idéal.
De plus le théorème \ref{theodistrib} devient :

\begin{theo}\label{theopdo} L'idéal $\cJ$ de $S_c(\T)$ est tel que pour tout $k=(k_t)_{t\in \R^+} \in \cJ$ et tout $m \in \N$ si :
 $$P(x,y)=\int_0^{+\infty} t^{-\mbox{dim}(M)} k_t(x,y)\dfrac{dt}{t^{1-m}} \ \mbox{ et }  \ \sigma(P) : (x,\xi)\in T^*M \mapsto \int_0^{+\infty} \hat{k_0}(x,t\xi)\dfrac{dt}{t^{1-m}}$$
 alors le noyau $P$ détermine un opérateur pseudodifférentiel (classique) d'ordre $-m$ sur $M$ de symbole principal $\sigma (P)$.\\
De plus tout opérateur pseudodifférentiel d'ordre négatif est de cette forme.
\end{theo}

\begin{rema}\label{suitePDO} Plus généralement, étant donné un groupoïde de Lie $G\rightrightarrows G^{(0)}$, on peut considérer le \emph{groupoïde adiabatique} : 
$$G_{ad}=\operatorname{DNC}(G,G^{(0)})_+=G\times \R^*_+\cup \frak{A}G\times\{0\} \rightrightarrows G^{(0)}\times \R_+$$
et définir de la même façon un calcul pseudodifférentiel sur $G$.

Les opérateurs pseudodifférentiels ainsi obtenus forment une algèbre de convolution. De plus :
\begin{itemize}
\item[\pt] Les opérateurs d'ordre $m<0$ sont dans $C^*(G)$,
\item[\pt] Les opérateurs d'ordre $m\leq 0$ sont des multiplicateurs (des opérateurs à droite et à gauche bornés) de $C^*(G)$.
\end{itemize}
On désigne par $\Psi^*(G)$ la fermeture de l'algèbre des opérateurs d'ordre $0$ dans la $C^*$-algèbre des multiplicateurs de $C^*(G)$. On obtient une suite exacte de $C^*$-algèbre :
$$\begin{CD}
   0 @>>> C^*(G) @>>> \Psi^*(G) @>\sigma>> C(\mathbb{S}\frak{A}^*G) @>>>0 \end{CD}$$ 

 où $\mathbb{S}\frak{A}^*G$ est le fibré en sphère du fibré dual de l'algébroïde de Lie $\frak{A}G$ et $\sigma$ l'application symbole principal.
\end{rema}

\smallskip Depuis les travaux précurseurs de Connes de nombreux auteurs ont prolongé  ces idées pour définir ou redéfinir des calculs pseudodifférentiels dans des contextes variés. Dans un cadre assez proche de celui qui nous concerne ici, citons les travaux de \textcite{Pon08,VanErpAS1, VanErpAS2}.

\section{Groupoïde tangent en géométrie sous-riemannienne et calcul pseudodifferentiel} 

Comme vu précédemment, le groupoïde tangent permet de définir le calcul pseudodifférentiel classique. La démarche d'Androulidakis, Mohsen, Yuncken consiste à appliquer une construction analogue sur un groupoïde tangent adapté au contexte de la géométrie sous-riemannienne.
Nous présentons dans cette partie la construction de ce groupoïde telle qu'elle apparaît dans \parencite{MoTanG, MoBlup}. Ce groupoïde servira de support, comme précédemment, à la construction d'un calcul pseudodifférentiel.

\smallskip Considérons à nouveau la variété $M$, la famille de champs de vecteurs $X_1,\dots,X_n$ satisfaisant la condition de Hörmander et $N\in \N$ tel que le nombre de commutateurs pour engendrer $T_xM$ est majoré par $N$ indépendant de $x$. On note $\cF$ le sous-module de $\cX(M)$ engendré par $X_1,\dots,X_n$ et $d_{CC}$ la distance de Carnot-Carathéodory correspondante (définition \ref{dcc}). Pour tout $X\in \cX(M)$ et $x\in M$, $\exp(X)x$ désigne le flot au temps 1 partant de $x$ du champ $X$.

Comme dans la section \ref{conetang}, $\frak{g}$ désigne l'algèbre de Lie libre engendrée par $\tilde{X_1},\dots, \tilde{X_n}$ avec pour relations que les suites de crochets de longueur strictement plus grande que $N$ s'annulent. Il s'agit d'une algèbre de Lie nilpotente de dimension finie. On note $G$ le groupe de Lie nilpotent simplement connexe qui intègre $\frak{g}$.

La formule de Baker-Campbell-Hausdorff offre une correspondance bijective entre les sous-algèbres de Lie de 
$\frak{g}$ et les sous-groupes simplement connexes de $G$. Pour $H$ un tel sous-groupe, on notera $\frak{h}\subset \frak{g}$ son algèbre de Lie et la bijection réciproque $\Xi(\frak{h})=H$.

\subsection{Groupoïde tangent en géométrie sous-riemannienne}\label{TangOmar} \label{Tau}

Une des idées essentielles dans le travail de Mohsen trouve sa source dans l'observation suivante. Fixons $x\in M$, les travaux de \parencite{Bel} permettent d'identifier un sous-groupe simplement connexe $G_x\subset G$ de codimension $\mbox{dim} M$ tel que $$\underset{t\rightarrow 0^+}{\mbox{lim}} (M,t^{-1}d_{CC},x)=(G/G_x, d_{G/G_x},G_x)$$
où $d_{G/G_x}$ est la métrique de Carnot-Carathéodory sur l'espace homogène $G/G_x$. 

Cette limite n'étant pas uniforme en $x$, elle ne suffit pas à définir l'\enquote{espace tangent} à la structure sous-riemannienne dans sa globalité. En effet, si $(x_n,t_n)_{n\in \N}$ est une suite de points de $M\times \R_+$ telle que $(x_n,t_n)$ converge vers $(x,0)$, quitte à prendre une sous-suite, la limite $\underset{n\rightarrow +\infty}{\mbox{lim}} (M,t_n^{-1}d_{CC},x_n)$ existe et est de la forme $(G/H,d_{G/H},H)$ où $H\subset G$ est un sous-groupe de Lie simplement connexe, mais en général $H\not= G_x$.

Enfin, moyennant l'identification des sous-groupes simplement connexes de $G$ avec leur algèbre de Lie, il apparaît que tous ces cônes tangents sont à chercher dans la variété Grassmannienne des sous-espaces de codimension $\mbox{dim} M$ de $\frak{g}$ ; nous noterons $\mathrm{Grass}(\frak g)$ cet espace.

\smallskip On peut noter, que cette difficulté n'apparaît pas dans le cas \emph{équirégulier} pour lequel le théorème de Helffer et Nourrigat était connu. On dit que la famille $\{X_1,\cdots,X_k\}$ est équirégulière si pour tout $k\in \N$, la dimension du sous-espace de $T_xM$ obtenu en évaluant $\cF^k$ en $x$ est localement constante. Dans ce cas, il n'y a pour chaque $x\in M$ qu'un seul cône tangent.

\paragraph*{\Gt L'espace des unités $\G^{(0)}$} Notons $\natural : \frak{g} \rightarrow \cX (M)$ l'application linéaire qui envoie chaque $\tilde{X_i}$ sur $X_i$, ainsi que les suites de crochets de longueur inférieure ou égale à $N$ sur les crochets correspondants de $\cX(M)$. Bien sûr, $\natural$ n'est pas un morphisme d'algèbre. 

\smallskip Pour $x\in M$, on considère $\natural_x= ev_x\circ \natural : \frak{g} \rightarrow T_x M$, où $ev_x$ est l'évaluation en $x$. Notons que la condition de Hörmander garantie la surjectivité de $\natural_x$ de sorte que son noyau est un élément de $\mathrm{Grass}(\frak g)$. Puisque $\natural$ n'est pas un morphisme d'algèbre de Lie, $\mbox{ker}(\natural_x)$ n'est pas une sous-algèbre de Lie de $\frak{g}$ en général.

On définit également la dilatation graduée $\alpha_t : \frak{g} \rightarrow \frak{g}$ donnée par : 
$\alpha_t(\tilde X_i)=t\tilde X_i$, $\alpha_t([\tilde X_i,\tilde X_j])=t^2[\tilde X_i,\tilde X_j]$...

Enfin, on considère l'inclusion 
$$\begin{array}{crcl} \Theta_0 : &  M\times \R_+^* & \rightarrow & \mathrm{Grass}(\frak{g})\times M \times \R_+ \\  & (x,t) & \mapsto & (\alpha_{\frac{1}{t}}(\mbox{ker}(\natural_x)),x,t) \end{array}$$
On définit l'espace :
$$\G^{(0)}:= M\times \R_+^* \underset{x\in M}{\sqcup} \cG_x^0\times \{ (x,0)\}$$
comme la fermeture de $M\times \R_+^*$ dans $\mathrm{Grass}(\frak{g})\times M \times \R_+$.

L'un des résultats principaux de \textcite{MoTanG} est le suivant 
\begin{theo} Si $(x_n,t_n)\in M\times \R_+^*$ converge dans $\G^{(0)}$ vers $(\frak{h},x,0)$, c'est à dire si $t_n \rightarrow 0$, $x_n \rightarrow x$ et $\alpha_{\frac{1}{t_n}}(\mbox{ker}(\natural_{x_n})) \rightarrow \frak{h}$ alors :
\begin{enumerate}
\item $\frak{h}$ est une sous-algèbre de Lie de $\frak{g}$.  

Notons $H=\Xi(\frak{h})$ le sous-groupe simplement connexe de $G$ d'algèbre de Lie $\frak{h}$.
\item Pour tout $g\in G$, l'algèbre de Lie de $gHg^{-1}$ est dans $\cG_x^0$.
\item On a  $$\underset{n\rightarrow +\infty}{\mbox{lim}} (M,t_n^{-1}d_{CC},x_n)=(G/H,d_{G/H},H)$$
\end{enumerate}
\end{theo}

De plus la projection $\G^{(0)} \rightarrow M\times \R_+$ est propre, de sorte que la construction précédente fournit tous les cônes tangents possibles.

\paragraph*{\Gt Lien avec l'ensemble de Helffer et Nourrigat $\cT_x$}

La méthode des orbites de Kirillov établit une correspondance entre les représentations irréductibles et unitaires d'un groupe de Lie et les orbites co-adjointes, c'est-à-dire, les orbites de l'action du groupe sur le dual de son algèbre de Lie. À l'aide de cette identification, on peut décrire le sous-ensemble $\cT_x$ des représentations  irréductibles et unitaires du groupe de Lie nilpotent $G$ comme un sous-ensemble de $\frak{g}^*$ des deux façons équivalentes suivantes.

\begin{defi} Le \emph{cône de Helffer et Nourrigat} en $x$, vu dans $\frak{g}^*$, est défini de façon équivalente par :
\begin{enumerate} 
\item L'élément $\xi$ de  $\frak{g}^*$ est dans $\cT_x$ si et seulement si il est limite d'une suite de la forme $(\xi_n\circ \natural_{x_n}\circ \alpha_{t_n})_{n\in \N}$ où $(t_n)_{n\in \N}$ est une suite d'éléments de $\R_+^*$ qui tend vers $0$, $(x_n,\xi_n)_{n\in \N}$ une suite d'éléments de $T^*M$ telle que la suite $(x_n)_{n\in \N}$ tend vers $x$.
\item  
$$\cT_x =\underset{\frak{h}\in \cG_x^0}{\bigcup} \frak{h}^\perp$$ où $\frak{h}^\perp=\{\xi \in \frak{g}^*  ; \ \xi(\frak{h})=0\}$.
\end{enumerate}
\end{defi}

La première définition provient des travaux de \textcite{HN79c}, la seconde se trouve dans les travaux de \textcite{AMY}. En outre, une définition alternative s'obtient par la proposition suivante.

\begin{prop}
Une représentation irréductible et unitaire $\pi$ du groupe de Lie nilpotent $G$  est dans $\cT_x$ si et seulement si $\pi$ est faiblement contenue dans la représentation induite $G\rightarrow U(L^2(G/H))$ pour un $H=\Xi(\frak{h})$ avec  $\frak{h} \in \cG_x^0$.
\end{prop}

\paragraph*{\Gt Le groupoïde $\G\rightrightarrows \G^{(0)}$} 

Dans la suite, on identifiera, les sous-groupes simplement connexes de $G$ et leur algèbre de Lie. Ainsi, pour tout $x\in M$, on identifie $\cG_x^0$ et $\Xi(\cG_x^0)$ de sorte que $\cG_x^0$ apparaît comme un ensemble de sous-groupes de Lie simplement connexes de $G$ et $$\G^{(0)}:= M\times \R_+^* \sqcup \{(H,x) \ ; x\in M, H\in\cG_x^0 \} \times\{0\}$$

\medskip Dans le cas du groupoïde tangent de Connes, on recolle $M\times M \times \R_+^*$ avec l'espace tangent $TM$ en $0$. En géométrie sous-riemannienne on procède de même, le rôle d'espace tangent étant cette fois joué par l'ensemble des cônes tangents. On pose donc : 
$$\G=M\times M \times \R_+^* \sqcup \{(gH,x) ; g\in G,\ x\in M,\ H\in\cG_x^0 \} \times \{0\} \rightrightarrows \G^{(0)}$$

La structure de groupoïde est la suivante :
\begin{itemize}
\item[\pt] Application unité :
$$u(x,t)=(x,x,t) \mbox{ si } t>0 \mbox{ et } u(H,x,0)=(H,x,0)$$
\item[\pt] Source et but :
$$\left\{\begin{array}{rcl} s(x,y,t)& = &(y,t)  \mbox{ si } t>0 \\ s(gH,x,0)& =& (H,x,0) \end{array} \right. \mbox{ et } \ \ \left\{\begin{array}{rcl} r(x,y,t)&=&(x,t)  \mbox{ si } t>0 \\ r(gH,x,0) &= &(gHg^{-1},x,0) \end{array} \right.$$
\item[\pt] Produit partiel :
$$\left\{\begin{array}{rcl} (x,y,t)\cdot (y,z,t) & = & (x,z,t)  \mbox{ si } t>0 \\  (hgHg^{-1},x,0)\cdot(gH,x,0) &= &( hgH,x,0) \end{array} \right.$$
\item[\pt] Inverse :
$$(x,y,t)^{-1}=(y,x,t) \mbox{ si } t>0 \mbox{ et } (gH,x,0)^{-1}=(g^{-1}(gHg^{-1}),x,0)=(Hg^{-1},x,0)$$
\end{itemize}

\smallskip On peut munir $\G$ d'une structure topologique naturelle de sorte que $\G$ est localement compact séparé à base dénombrable et un groupoïde $C^{\infty,0}$ d'unités $\G^{(0)}$. Plus encore, on peut construire un plongement topologique fermé de $\G$ dans $\R^n$ pour un $n$ assez grand.

\medskip La topologie obtenue sur ce groupoïde est la suivante :

Une suite $(y_n,x_n,t_n)_{n\in \N}$ converge vers $(gH,x,0)$ si et seulement si :
\begin{enumerate} \item $(x_n)_{n\in \N}$ et $(y_n)_{n\in \N}$ convergent vers $x$ et $(t_n)_{n\in \N}$ converge vers $0$,
\item la suite $(x_n,t_n,0)_{n\in \N}$ converge vers $(H,x,0)$ dans $\G^{(0)}$, 
\item Il existe $v\in gH$ et $(v_n)_{n\in \N}$ une suite de $G$ qui converge vers $v$ tels que, moyennant l'identification de $G$ avec son image dans $\cX(M)$ (via l'exponentielle qui identifie $G$ avec $\frak{g}$ puis   $\natural : \frak{g} \rightarrow \cX (M)$) on a $\exp(\alpha_{t_n}(v_n)).x_n=y_n$. \end{enumerate}

\begin{rema} Dans leur article \textcite{AMY} considèrent le feuilletage singulier $\frak{a} \cF := t\tilde{\cF_1}+t^2\tilde{\cF_1} +\cdots + t^N\tilde{\cF_N}$ sur $M\times \R$ où $\tilde{\cF_k}$ désigne les champs de vecteurs sur $M\times \R$ de la forme $\tilde{X}(x,t)=X(x)$ où $X\in \cF_k$. Ils construisent le groupoïde d'holonomie de ce feuilletage \parencite{AndrSk1} : 
$$\cH(\frak{a} \cF)=M\times M\times \R^* \underset{x\in G_x}{\sqcup} G_x\times \{0\} \rightrightarrows M\times \R$$
où $G_x$ est le groupe de Lie simplement connexe qui intègre l'algèbre de Lie $\frak{g}_x=\frak{a} \cF/I_{(x,0)} \frak{a} \cF$. Ils se concentrent ensuite sur l'adhérence de $M\times M\times \R_+^*$ dans ce groupoïde (singulier). Ce point de vue permet de travailler dans des \enquote{cartes} adaptées.
\end{rema}

\subsection{Calcul pseudodifférentiel et preuve de la conjecture de Helffer et Nourrigat}

\paragraph*{ \Gt Le symbole est correctement défini}
La première étape consiste à montrer que le symbole principal est correctement défini. Le théorème suivant utilise de façon cruciale  le fait que la $C^*$-algèbre de $\G$ est une $C_0(\R_+)$-algèbre ce qui permet d'obtenir que si $a\in C^*(\G)$ : $$\underset{t\rightarrow 0^+}{\mathrm{lim\ sup}}\|a_t\|_{K(L^2(M)}=\underset{x\in M}{\mathrm{sup}}\ \underset{\pi\in \cT ^*_x}{\mathrm{sup}} \|\pi(a_{x,0})\|_{L^2(\pi)}$$

On a :
\begin{theo} $[$\cite[Théorème 1.15]{AMY} $]$\\
Soit $D$ un opérateur différentiel et $x\in M$, pour tout $\pi \in \cT_x$, $\sigma(D,x,\pi)$ est bien défini.
\end{theo}

\paragraph*{\Gt Définitions des opérateurs pseudodifférentiels}
On procède ensuite à la définition d'un calcul pseudodifférentiel, à l'aide du groupoïde tangent en géométrie sous-riemannienne $\G$ de façon analogue à ce que l'on a fait en \ref{pdoclas}. Voici une idée de la construction.

Notons $\cG^0=\underset{x\in M}{\sqcup}\{(H,x)\ ; x\in M,\ H\in \cG_x^0\}$. On a (localement) une \enquote{submersion $C^{\infty,0}$} surjective:
$$\begin{array}{cccc} \Upsilon : & \cG^0 \times \R_+ \times \frak{g} & \rightarrow & \G \\  & (H,x,t,v) & \mapsto & (\exp(\alpha_t(\natural(v)))\cdot x,x,t) \mbox{ si } t\not=0 \\
& (H,x,0,v) & \mapsto & (\Xi(v)H,x,0) \end{array} $$

où $\Xi$ est l'application exponentielle de $\frak{g}$ dans $G$. 

Le membre de gauche est un fibré vectoriel sur $\cG^0 \times \R_+$ pour lequel on peut construire un $\cJ(\cG^0 \times \R_+ \times \frak{g},\cG^0 \times \R_+)$ (bien que $\cG^0$ ne soit pas une variété, mais un fermé dans une variété).

\smallskip Une fonction $f\in C_c^{\infty} ( \G)$ est dans $\cJ(\cF^{\bullet})$ si et seulement si $f\circ \Upsilon$ est dans $\cJ(\cG^0 \times \R_+ \times \frak{g},\cG^0 \times \R_+)$.

\begin{defi}  Pour $m\in \N$, on définit $\Psi^{-m}(\cF^{\bullet})$ comme l'ensemble des opérateurs de noyaux de la forme : 
 $$P(x,y)=\int_0^{+\infty} w(x,y,t) k_t(x,y)\dfrac{dt}{t^{1-m}}$$ 
 ou $k=(k_t)_{t\in \R^+} \in \cJ(\cF^{\bullet})$ et $w$ est un terme correctif\footnote{La fonction $w$ n'apparaît pas si l'on utilise des densités.} .  
\end{defi}

\medskip Cela permet d'étendre l'algèbre $\N$-filtrée des opérateurs différentiels $DO^m(\cF^{\bullet})$ en une algèbre $\Z$-filtrée $\Psi^{m}(\cF^{\bullet})$. 

Avec les notations de la définition ci-dessus, le symbole principal de l'opérateur $P$ est donné par $$\begin{aligned} \sigma(P,x,\pi)(\xi)= \pi( \int_0^{+\infty} (\alpha_t k_0)(x,\xi)\dfrac{dt}{t^{1-m}}) \end{aligned}$$ de plus le symbole d'un produit est le produit des symboles.

\smallskip \paragraph*{ \Gt L'algèbre $\Psi(\cF^{\bullet})$ a les propriétés \enquote{attendues} d'un calcul pseudodifférentiel}

\begin{theo}\label{theo1} $[$\cite[Théorème 3.31]{AMY} $]$\\
Soit $P\in \Psi^{k}(\cF^{\bullet})$ :
\begin{itemize}\item[\pt] Pour $k\leq 0$, $P$ s'étend en un opérateur borné de $L^2(M)$.
\item[\pt] Pour $k< 0$, $P$ s'étend en un opérateur compact de $L^2(M)$.
\end{itemize}
\end{theo}

En particulier, en étendant les espaces de Sobolev de la définition \ref{SobetOp} à $\Z$ par dualité : pour tout $k\in \Z$, $P\in \Psi^{k}(\cF^{\bullet})$ et $s\in \Z$, $P$ s'étend en un opérateur $P:
\tilde{H}^{s+k} \rightarrow \tilde{H}^{s}$.

\medskip Ce théorème permet de montrer que $\Psi(\cF^{\bullet})$ satisfait la \enquote{propriété de \emph{régularité}} :
\begin{itemize} \item[\pt] Les éléments de $\Psi^{-\infty}(\cF^{\bullet}):=\underset{m\in \Z}{\bigcap} \Psi^m(\cF^{\bullet})$ sont des opérateurs régularisants.
\end{itemize}

\bigskip Il découle de la définition des opérateurs pseudodifférentiels que $\Psi(\cF^{\bullet})$ satisfait la \enquote{propriété de \emph{pseudo-localité}} :
\begin{itemize} \item[\pt] Les noyaux de Schwartz des éléments de $\Psi(\cF^{\bullet})$ sont lisses en dehors de la diagonale.
\end{itemize}

\begin{theo}\label{theo3} $[$\cite[Proposition 3.32]{AMY} $]$\\
Pout tout $k\in \Z$, il existe $P_k\in  \Psi^{k}(\cF^{\bullet})$ et $Q_k\in  \Psi^{-k}(\cF^{\bullet})$ tels que $P_kQ_k-Id$ et $Q_kP_k-Id$ sont dans $\underset{m\in \Z}{\cap} \Psi^m(\cF^{\bullet})$.
\end{theo}

\smallskip Le théorème suivant utilise le théorème précédent pour se ramener à $\Psi^{0}(\cF^{\bullet})$, puis se déduit d'une suite exacte courte de $C^*$-algèbres analogue à la suite exacte des opérateurs pseudodifférentiels classiques présentée dans la remarque \ref{suitePDO}.
\begin{theo}\label{theo2} $[$\cite[Théorème 3.38]{AMY} $]$\\
Soit $P\in \Psi^{k}(\cF^{\bullet})$, les assertions suivantes son équivalentes :
\begin{enumerate}
\item Pour tout $x\in M$ et $\pi \in \cT_x\setminus{1_G}$, $\sigma(P,x,\pi)=0$.
\item Pour tout (ou un) $s\in \Z$ l'opérateur $P: \tilde{H}^{s+k} \rightarrow \tilde{H}^{s}$ est compact.
\end{enumerate}
\end{theo}

\begin{rema} Il est important de noter que si $P\in \Psi^{k}(\cF^{\bullet})$ et que pour tout $x\in M$ et $\pi \in \cT_x\setminus{1_G}$, $\sigma(P,x,\pi)=0$, on ne peut pas en déduire pour autant que $P\in \Psi^{k-1}(\cF^{\bullet})$. Ce point est une différence majeure avec le calcul pseudodifférentiel classique et rend particulièrement délicate la quête de paramétrix.  \end{rema}

 \smallskip Le théorème \ref{theo3} permet de réduire l'étude de l'existence de paramétrix aux éléments de $\Psi^{0}(\cF^{\bullet})$ ; éléments qui s'étendent d'après le théorème \ref{theo1} en des opérateurs bornés de $L^2(M)$ dans $L^2(M)$.

\begin{theo}\label{theo4} $[$\cite[Théorèmes 3.43]{AMY} $]$\\
Soit  $P\in \Psi^{0}(\cF^{\bullet})$. Les assertions suivantes sont équivalentes :
\begin{enumerate}
\item Pour tout $x\in M$ et $\pi\in \cT_x\setminus\{1_G\}$, $\sigma(P,x,\pi)$ est injective.
\item Il existe $Q\in \Psi^{0}(\cF^{\bullet})$ tel que pour tout $x\in M$ et $\pi\in \cT_x\setminus\{1_G\}$, $\sigma(QP,x,\pi)=Id$ 
\item Il existe un opérateur $Q:C^{-\infty}(M)\rightarrow C^{-\infty}(M)$ qui définit pour tout $s\in \Z$ un opérateur borné 
$Q:\tilde{H}^s\rightarrow \tilde{H}^{s}$ et tel que $QP-Id$ est régularisant.
\end{enumerate}
\end{theo}

\begin{proof}[Quelques mots sur l'implication de (2) à (3)]
Si un opérateur $P\in  \Psi^{0}(\cF^{\bullet})$ est tel que pour tout $x\in M$ et $\pi\in \cT_x\setminus\{1_G\}$, $\sigma(P,x,\pi)=Id$ alors, d'après le théorème \ref{theo2}, pour tout $s\in \Z$ l'opérateur $Id-P:\tilde{H}^s\rightarrow \tilde{H}^s$ est compact. Donc, pour tout $s\in \Z$, $P:\tilde{H}^s\rightarrow \tilde{H}^s$ est un opérateur de Fredholm d'indice de Fredholm nul. Quitte à lui ajouter un régularisant $R\in C^{\infty}(M\times M)$ il est inversible comme opérateur de $L^2(M)$ dans $L^2(M)$. Notons $Q=(P+R)^{-1} :L^2(M) \rightarrow L^2(M)$.

Regardons $P+R:\tilde{H}^s\rightarrow \tilde{H}^s$. Lorsque $s>0$, puisque $\tilde{H}^s \subset L^2(M)$, l'opérateur $P+R$ est injectif, de Fredholm et d'indice nul, il est donc inversible. Lorsque $s<0$, puisque $L^2(M)\subset  \tilde{H}^s$, l'opérateur $P+R$ est surjectif, d'image dense puisqu'elle contient 
$L^2(M)$ et donc inversible à nouveau. Donc pour tout $s\in \Z$, $Q$ s'étend de $\tilde{H}^s\rightarrow \tilde{H}^s$.

En utilisant que les $C^*$-algèbres sont stables par calcul fonctionnel holomorphe, on obtient que pour tout $s\in \Z$, $Q=(P+R)^{-1}$ est dans la fermeture de $\Psi^{0}(\cF^{\bullet})$ dans $\cL(\tilde{H}^s)$. \end{proof}

\smallskip On en déduit que les éléments que $\Psi(\cF^{\bullet})$ admettent des \enquote{\emph{paramétrix} }:
\begin{itemize} \item[\pt]
Si le symbole principal de $P \in \Psi^k(\cF^{\bullet})$ est inversible, il existe un opérateur $Q: C^{-\infty}(M)\rightarrow C^{-\infty}(M)$ qui définit pour tout $s\in \Z$ un opérateur borné $Q:\tilde{H}^s\rightarrow \tilde{H}^{s+k}$ et tel que $QP-Id$ et $QP-Id$ sont régularisants.
\end{itemize}

\begin{rema} Contrairement au calcul pseudodifférentiel classique, les paramétrix peuvent ne pas être des éléments de $\Psi(\cF^{\bullet})$. 
\end{rema}

\smallskip On est alors en mesure de montrer que l'inversibilité du symbole principale d'un opérateur différentiel implique l'hypoellipticité. En effet, pour $k\in \N$, si le symbole principal de $D \in DO^k(\cF^{\bullet})$ est inversible, considérons un paramétrix $Q$ tel que $QD-Id \in \Psi^{-\infty}(\cF^{\bullet})$.
 
 Si $u \in C^{-\infty}(M)$, on écrit $u=(I-QD)u+QDu$ et on a :
 \begin{itemize}
 \item[\pt] $(Id-QD)u$ est lisse car $Id-QD$ est régularisant.
 \item[\pt] Si $Du$ est dans $C^{\infty}(M)=\underset{s\in \Z}{\bigcap} \tilde{H}^s$ alors $QDu$ aussi, puisque pour tout $s\in \Z$, $Q$ induit un opérateur $Q:\tilde{H}^s\rightarrow \tilde{H}^{s+k}$.
 \end{itemize}

Finalement si $f$ est lisse et $Du=f$ alors $u$ est lisse.

\paragraph*{\ding{71} Un exemple} La variété considérée est $\R^2$ avec $X_1=\partial_x$ et $X_2=x\partial_y$\footnote{Pour être vraiment dans une variété compacte, on peut considérer $M=S^1\times S^1$ avec $X_1=\partial_x$ et $X_2=\sin(x)\partial_y$.}. \\
On a $[X_1,X_2]=\partial_y$ donc $N=3$ et l'algèbre de Lie $\frak{g}$ est engendrée (comme espace vectoriel) par trois éléments $\tilde{X_1}, \tilde{X_2}$ et $\tilde{X_3}=[\tilde{X_1},\tilde{X_2}]$ qui sont tels que $[\tilde{X_1},\tilde{X_3}] =[\tilde{X_2},\tilde{X_3}]=0$. En particulier $Grass(\frak{g})=\mathbb{P}^1(\frak{g})$. \\ On reconnait l'algèbre de Lie du groupe de Heisenberg, donc $G=\mathbb{H}$. 

Rappelons que l'on a les identifications suivantes : 

$G=\{\left(\begin{smallmatrix} 1 & a_1 & a_3 \\ 0 & 1 & a_2 \\ 0 & 0 & 1\end{smallmatrix} \right); a_i \in \R \}$, 
$\frak{g}=\{\left(\begin{smallmatrix} 0 & v_1 & v_3 \\ 0 & 0 & v_2 \\ 0 & 0 & 0\end{smallmatrix} \right) ; v_i \in \R\}$, $\frak{g}^*=\{A^t \ ; A\in \frak{g}\}$ avec pour $\xi \in \frak{g}^*$ et $v\in \frak{g}$ : $\langle \xi,v\rangle=\mathrm{Tr} (\xi v)$ .

Pour $g=\left(\begin{smallmatrix}1 & a_1 & a_3 \\ 0 & 1 & a_2 \\ 0 & 0 & 1\end{smallmatrix} \right)\in G$ et 
$\xi=\left(\begin{smallmatrix} 0 & 0 & 0 \\ f_1 & 0 & 0 \\ f_3 & f_2 & 0\end{smallmatrix} \right) \in \frak{g}^*$, l'action co-adjointe de $g$ sur $\xi$ est donnée par $g {\scriptscriptstyle \bullet} \xi=\left(\begin{smallmatrix} 0 & 0 & 0 \\ f_1 +a_2f_3& 0 & 0 \\ f_3 & f_2-a_1f_3 & 0\end{smallmatrix} \right)$. On a donc des orbites de dimension $2$ lorsque $f_3\not=0$ et de dimension $0$ sinon.

\smallskip Les représentations unitaires irréductibles du groupe de Heisenberg $G$ sont de deux types :
\begin{itemize}
\item[\pt] Pour $\varepsilon \in \{-1,1\}$, $\pi_{\varepsilon}:G\rightarrow U(L^2(\R))$ est engendrée au niveau de son algèbre de Lie par  :
$$\tilde{X_1}\mapsto \varepsilon\partial_t,\ \tilde{X_2}\mapsto  it,\ \tilde{X_3}\mapsto i\varepsilon$$ et correspondent aux orbites de dimension 2 de l'action co-adjointe de $G$ sur $\frak{g}^*$.
\item[\pt] Pour $(a,b)\in \R^2\setminus\{(0,0)\}$, $\pi_{a,b} : G \rightarrow U(\C)$ est engendrée au niveau de son algèbre de Lie par  :
$$\tilde{X_1}\mapsto ia,\ \tilde{X_2}\mapsto ib,\ \tilde{X_3}\mapsto 0$$ et correspondent aux orbites de  dimension 0 de l'action co-adjointe de $G$ sur $\frak{g}^*$.
\end{itemize}

\medskip Pour $(x,y)\in \R^2$, $\natural_{(x,y)} :  \left(\begin{smallmatrix} 0 & v_1 & v_3 \\ 0 & 0 & v_2 \\ 0 & 0 & 0\end{smallmatrix} \right)\in \frak{g} \mapsto v_1\partial_x+(v_2x+v_3)\partial_y\in T_{(x,y)}\R^2$.

\smallskip Pour $t\in \R_+^*$, $(x,y)\in \R^2$, $\alpha_{\frac{1}{t}}ker(\natural_{(x,y)})= \langle  \left(\begin{smallmatrix} 0 & 0 & -x \\ 0 & 0 & t \\ 0 & 0 & 0\end{smallmatrix} \right) \rangle$. 

\smallskip On s'intéresse à la limite de  $\alpha_{\frac{1}{t_n}}ker(\natural_{(x_n,y)})$ lorsque $t_n$ tend vers $0$ et $x_n$ tend vers $x$.
\begin{itemize}
\item[\pt] Pour $x\not=0$ cette limite est toujours $ \langle \tilde{X_3} \rangle=\langle\left(\begin{smallmatrix} 0 & 0 & 1 \\ 0 & 0 & 0 \\ 0 & 0 & 0\end{smallmatrix} \right)\rangle$.
Donc $\cG^0_{x,y}=\langle \tilde{X_3} \rangle$ et $\cT_{x,y}=\langle \tilde{X_3} \rangle^\perp =\{\left(\begin{smallmatrix} 0 & 0 & 0 \\ f_1 & 0 & 0 \\ 0 & f_2 & 0\end{smallmatrix} \right) \ ; f_i\in \R\} \subset \frak{g}^*$.

Les représentations à considérer sont donc celles du type $\pi_{a,b}$ pour $(a,b)\in  \R^2$.

\item[\pt] Pour $x=0$, en notant que $\langle  \left(\begin{smallmatrix} 0 & 0 & -x_n \\ 0 & 0 & t_n \\ 0 & 0 & 0\end{smallmatrix} \right) \rangle = \langle  \left(\begin{smallmatrix} 0 & 0 & -\frac{x_n}{t_n} \\ 0 & 0 & 1 \\ 0 & 0 & 0\end{smallmatrix} \right) \rangle$, on obtient alors que  $\cG^0_{x,y}=\langle \tilde{X_3}   \rangle \underset{\lambda\in \R}{\sqcup} \langle \tilde{X_2}+\lambda \tilde{X_3}   \rangle$.

Donc $\cT_{x,y}=\{\left(\begin{smallmatrix} 0 & 0 & 0 \\ f_1 & 0 & 0 \\ 0 & f_2 & 0\end{smallmatrix} \right) \ ; f_i\in \R\} \underset{\lambda\in \R}{\sqcup} \{\left(\begin{smallmatrix} 0 & 0 & 0 \\ f_1 & 0 & 0 \\ f_3 & -\lambda f_3 & 0\end{smallmatrix} \right) \ ; f_1,f_3\in \R\}=\frak{g}^*$ .

Les représentations à considérer sont donc les $\pi_{a,b}$ pour $(a,b)\in  \R^2$ ainsi que les $\pi_\varepsilon$ pour $\varepsilon \in \{-1,1\}$.
\end{itemize}

\medskip On se donne une fonction $\ell\in C_c^{\infty}(\R^2)$ et on considère l'opérateur différentiel sur $\R^2$ : $$D_\ell=\partial_x^2+x^2\partial_y^2+\ell(x,y)i\partial_y$$ On a $D_\ell=P(X_1,X_2)$ où $P(X,Y)=X^2+Y^2+\ell(x,y)i(XY-YX) $ et $D$ est d'ordre de Hörmander $2$.

On a, pour $(a,b)\in \R^2\setminus \{(0,0)\}$ et $\varepsilon \in \{-1,1\}$ :
\begin{itemize} \item[\pt] $\sigma(D_\ell,(x,y),\pi_{a,b})=((ia)^2+(ib)^2)Id=-(a^2+b^2)Id$. 

On retrouve le symbole classique.
\item[\pt] $\sigma(D_\ell,(0,y),\pi_\varepsilon)= \partial_t^2+((it)^2-\ell(0,y)\varepsilon) Id_{L^2(\R)}=(\partial_t^2-t^2)-\ell(0,y)\varepsilon Id_{L^2(\R)}$
\end{itemize}

\smallskip On conclut que l'opérateur $D_\ell$ est hypoelliptique maximal si et seulement si \mbox{$\ell(0,y) \notin \varepsilon\cdot \mbox{spec}(\partial_t^2-t^2)=\{2n+1,\ n\in \Z\}$}.






\printshorthands 

\printbibliography

\end{document}